\documentclass[letterpaper,10pt,conference]{ieeeconf}
\IEEEoverridecommandlockouts
\usepackage[T1]{fontenc}
\usepackage{amsmath}
\usepackage{newtxmath}
\usepackage{bm}
\usepackage{cite}
\usepackage{xurl}
\usepackage{graphicx}
\usepackage{array}
\usepackage{booktabs}
\usepackage{tikz}
\usetikzlibrary{arrows.meta,positioning,fit,backgrounds}
\makeatletter\let\NAT@parse\undefined\makeatother
\usepackage[hidelinks]{hyperref}
\hypersetup{pdftitle={Optimal Control of Dynamic Gas Networks Using Gaussian Processes with Chance-Constraint Reformulation},pdfauthor={Taiwo A. Adebiyi, Robert Ferrando, Manuel Garcia, Ruda Zhang, Saif R. Kazi, Kaarthik Sundar}}
\newcolumntype{L}[1]{>{\raggedright\arraybackslash}p{#1}}
\title{\LARGE\bfseries Optimal Control of Dynamic Gas Networks Using Gaussian Processes with Chance-Constraint Reformulation}
\author{Taiwo A.~Adebiyi$^{1,2}$, Robert Ferrando$^{3, *}$, Manuel Garcia$^{3}$, Ruda Zhang$^{1}$,
Saif R.~Kazi$^{2}$ and Kaarthik Sundar$^{3}$%
\thanks{$^{1}$Department of Civil and Environmental Engineering, University of Houston, Houston, TX 77204,
USA {\tt\footnotesize  \{taadebiyi2,\allowbreak rudaz\}@uh.edu}}%
\thanks{$^{2}$Applied Mathematics \& Plasma Physics, Los Alamos National Laboratory, Los Alamos, NM 87545,
USA {\tt\footnotesize skazi@lanl.gov}}%
\thanks{$^{3}$Computational Intelligence \& Modeling, Los Alamos National Laboratory, Los Alamos, NM 87545, USA {\tt\footnotesize \{rferrando,\allowbreak mjgarcia,\allowbreak kaarthik\}@lanl.gov}. $^{*}$Corresponding author. \mbox{LA-UR-26-28461}.}%
}

\begin{document}
\raggedbottom

\maketitle
\suppressfloats[t]
\thispagestyle{empty}
\pagestyle{empty}

\begin{abstract}
Natural gas networks must accommodate variable demand and supply while maintaining pressure across pipelines whose stored gas couples operating decisions through time. Probabilistic surrogates can simplify compressor scheduling, but their predictive uncertainty must remain reliable where optimization uses it. We combine direct learning of the minimum pressure margin with assessment at optimizer-selected schedules. Using a warped Gaussian process, we express the posterior pressure requirement as one inequality in its predictive mean and standard deviation. We then parameterize block and Fourier schedules so scalar predictions and analytic derivatives replace trajectory states and time integration in whole-day optimization. In steady state, a surrogate that passes held-out accuracy and coverage checks nevertheless falls below the minimum pressure on $71\%$ of its own dispatches at a nominal $1\%$ tolerance, as the optimizer is drawn to schedules where the margin is overestimated. This diagnostic motivates an enrichment loop that audits selected dispatches against the physics and updates the surrogate by exact conditioning, with fitted quantities held fixed. The loop restores steady-state minimum-pressure risk control at every studied tolerance, with conditioning about four orders of magnitude faster than refitting. In transient operation, we complement this enrichment with empirical threshold calibration to address residual hourly optimism, bringing hourly violation rates closer to their declared tolerances with little compressor work overhead.
\end{abstract}

\section{Introduction}
\label{sec:introduction}

Natural gas transmission presents a modeling tradeoff: the gas stored in pipelines, or linepack, provides operating flexibility, but exploiting it requires operating decisions, such as compressor schedules, whose effects extend across the network and through time. Established methods solve the coupled pressure and flow equations in steady and transient operation to evaluate decisions and optimize operation \cite{srinivasan2023,zlotnik2015}. The challenge is not the absence of physical models, but their repeated use across changing operating conditions and alternative compressor schedules. Learned surrogates offer a complementary approach, constructing reusable approximations of pipe relations \cite{arrigo2022} and network dynamics \cite{wang2025} from simulation or measurement data. Their computational usefulness, however, does not settle how reliably operating limits can be enforced through their predictions.

\begin{figure}[t]
\centering
\includegraphics{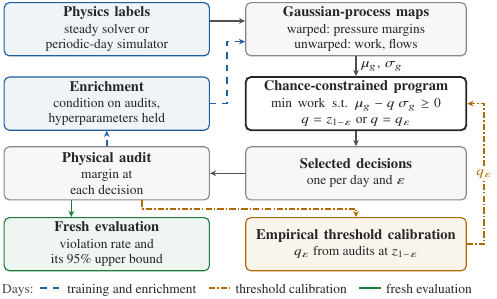}
\caption{Learned-margin workflow. The optimizer selects a least-work decision satisfying $\mu_g\ge q\sigma_g$ for each day and tolerance $\varepsilon$, then the physics audits it. Here $\mu_g$ and $\sigma_g$ are the transformed-margin posterior mean and standard deviation, with $q=z_{1-\varepsilon}$ by default. Enrichment conditions on audits with fitted parameters fixed (dashed loop). Transient calibration may replace $q$ by $q_\varepsilon$; fresh days evaluate violation rates.}
\label{fig:framework}
\end{figure}

The natural gas community has addressed uncertain operation through robust scheduling \cite{roald2020}, chance-constrained steady-state control \cite{dvorkin2022}, and stochastic network-flow representations \cite{kazi2024sfv}. These methods protect against withdrawal ranges or propagate uncertain inputs through physical equations. A surrogate introduces different uncertainty: how well its prediction represents the physical response, even with withdrawals and supply pressure prescribed. When that prediction determines admissibility, what the surrogate learns becomes part of the formulation.

Existing gas surrogates approximate pipe-flow relations \cite{arrigo2022}, nodal pressures and flows \cite{jiang2026}, or transient state-space dynamics \cite{wang2025}. These representations support
many downstream tasks, but a network-wide minimum-pressure requirement concerns an extremum: the
minimum \emph{pressure margin} across nodes and, during transient operation, across time. One could learn the
pressures and then derive this minimum, but its probability law depends on their high-dimensional
joint distribution. Learning the margin directly preserves the violation event in one scalar, regardless of the limiting node or instant, consistent with limit-state approximation in reliability-based design \cite{dubourg2011}.

We use a Gaussian process (GP) \cite{rasmussen2006} to pose the pressure requirement as a chance constraint with a declared violation tolerance. Chance constraints are established in control and network operation \cite{mesbah2016}, including GP-constrained optimization \cite{wiebe2022} and stability-constrained optimal power flow \cite{divito2026}. We develop a warped GP using a strictly increasing, zero-preserving transformation to accommodate margin asymmetry while retaining the known exact deterministic equivalent of the posterior chance constraint \cite{snelson2003,nemirovski2006}.

\begin{figure}[!t]
\centering
\includegraphics[width=\columnwidth]{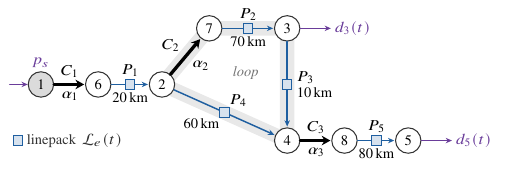}
\caption{The eight-node network used in both studies. Gas enters at slack node~1 and leaves at nodes~3 and~5. Heavy arrows denote compressors $C_1$--$C_3$; arrows show reference orientations, and gray shading marks loop 2--7--3--4. Squares mark pipe linepack $\mathcal L_e(t)$; purple identifies prescribed supply pressure $p_s$ and withdrawals $d_3(t)$ and $d_5(t)$.}
\label{fig:network}
\end{figure}

Although GPs offer principled uncertainty \cite{rasmussen2006}, neither their construction nor warping guarantees calibration, the agreement between predicted probabilities and observed frequencies \cite{gneiting2007}. Accurate predictions can understate uncertainty \cite{kuleshov2018,yadav2026}, and average interval coverage can conceal lower-tail errors \cite{allen2025}. We therefore assess accuracy, bulk coverage, and tail coverage among low-margin cases. This multi-axis evaluation of uncertainty has practical implications for compressor scheduling: understated uncertainty can admit pressure shortfalls, while overstated uncertainty can require unnecessary compression. Predictions also guide where the model is evaluated, so errors at optimizer-selected schedules need not resemble those on a random validation design \cite{fannjiang2022}. The question is whether the learned probability remains informative at selected decisions.

We address this question by coupling direct pressure-margin learning with physical audits of optimizer-selected schedules (Fig.~\ref{fig:framework}). We condition the GP on margins returned by the reference solver or simulator,
holding its fitted hyperparameters, warp parameters, and preprocessing fixed. In the steady-state diagnostic, this update improves
accuracy at selected decisions and restores minimum-pressure risk control. Daily transient operation provides the practical setting: block levels or Fourier-generated hourly setpoints parameterize schedules, and the learned whole-day margin replaces trajectory states in the optimizer. Exact conditioning and empirical threshold calibration bring audited hourly violation rates closer to declared tolerances on fresh days, with little compressor work overhead.

Sections~\ref{sec:problem} and~\ref{sec:method} develop the control problem and learned-margin optimization;
Sections~\ref{sec:steady} and~\ref{sec:transient} report the steady and transient studies, and Section~\ref{sec:conclusion} concludes.
\section{Gas-network control problem}
\label{sec:problem}

\subsection{Network physics and dispatch}
\label{sec:physics}

We consider a pipeline network $\mathcal G=(\mathcal N,\mathcal E)$ (Fig.~\ref{fig:network}), with nodes $\mathcal N$ and edges $\mathcal E=\mathcal P\cup\mathcal C$, comprising pipes $\mathcal P$ and compressors $\mathcal C$. Positive flow follows each edge's reference orientation $e=(i,j)$ from node $i$ to node $j$. For horizontal pipes in the low-velocity transmission regime, we neglect gravity and convective acceleration while retaining the temporal change of mass flux \cite{gyrya2019}. With the isothermal ideal-gas relation $p_e=a^2\rho_e$, where $a$ is the isothermal wave speed, the mass and momentum balances are
\begin{subequations}\label{eq:gas_dynamics}
\begin{align}
\partial_t\rho_e+\partial_\ell\phi_e&=0,\label{eq:mass_dynamics}\\
\partial_t\phi_e+\partial_\ell p_e&=-\frac{\lambda_e}{2D_e}\frac{\phi_e|\phi_e|}{\rho_e},\label{eq:momentum_dynamics}
\end{align}
\end{subequations}
where $\rho_e$, $\phi_e$, and $p_e$ are density, mass flux per unit area, and pressure; $t$ is time, $\ell\in[0,L_e]$ is distance along a pipe of length $L_e$, and $\lambda_e$ and $D_e$ are its friction factor and diameter. Integrating \eqref{eq:mass_dynamics} over a pipe of cross-sectional area $A_e$ gives its stored gas, or \emph{linepack}, and its rate of change,
\begin{equation}
\mathcal L_e(t)=A_e\int_0^{L_e}\rho_e(\ell,t)\,d\ell,
\qquad \dot{\mathcal L}_e=f_e^{\mathrm{in}}-f_e^{\mathrm{out}},
\label{eq:linepack}
\end{equation}
where $f_e^{\mathrm{in}}$ and $f_e^{\mathrm{out}}$ are inlet and outlet mass flows. Dependence on earlier injection and compression motivates optimizing transient rather than independent steady-state operation \cite{zlotnik2015}.

In steady state, \eqref{eq:mass_dynamics} makes $f_e=A_e\phi_e$ constant along each pipe. Squared pressures $\Pi_j=p_j^2$ express the integrated momentum balance as a nodal difference \cite{srinivasan2023}. Compressor $c=(i,j)$ has ratio $\alpha_c=p_j/p_i$. Let $s$ denote the slack node at prescribed pressure $p_s$, and $\mathcal N_0=\mathcal N\setminus\{s\}$. For net withdrawal $d_j$ and entering and leaving edge sets $\delta^{\mathrm{in}}(j)$ and $\delta^{\mathrm{out}}(j)$, the network equations are
\begin{subequations}\label{eq:steady_physics}
\begin{align}
\Pi_i-\Pi_j&=\kappa_e f_e|f_e|, &&e\in\mathcal P,\label{eq:pipe_law}\\
\Pi_j&=\alpha_c^2\Pi_i, &&c\in\mathcal C,\label{eq:compressor_relation}\\
d_j&=\textstyle\sum_{e\in\delta^{\mathrm{in}}(j)}f_e-\sum_{e\in\delta^{\mathrm{out}}(j)}f_e, &&j\in\mathcal N_0,\label{eq:nodal_balance}\\
\Pi_s&=p_s^2,\label{eq:supply_pressure}
\end{align}
\end{subequations}
where $\kappa_e=a^2\lambda_eL_e/(D_eA_e^2)$ is the pipe resistance. In transient operation, the nodes store no gas, so \eqref{eq:nodal_balance} balances the instantaneous pipe-end and compressor flows, and \eqref{eq:compressor_relation} links inlet and outlet pressures at each instant.

In this study, we collect the prescribed inputs in $\bm\xi$, comprising withdrawals and slack pressure, constant in steady state, or withdrawal profiles and constant slack pressure over a transient day $[0,T]$. We minimize compression work over steady ratios $\bm\alpha$ or coefficients $\bm\theta$ defining transient schedules $\alpha_c(t;\bm\theta_c)$, where $\bm\theta_c$ collects compressor $c$'s coefficients. For isentropic ideal-gas compression with unit isentropic efficiency, the work rate at mass flow $f_c$ is
\begin{equation}
\dot W_c(\alpha_c,f_c)=\eta f_c\bigl(\alpha_c^{(\gamma-1)/\gamma}-1\bigr),
\label{eq:compression_power}
\end{equation}
where $\eta$ is the specific-work constant and $\gamma$ the heat-capacity ratio \cite{kazi2024sfv,roald2020}. We therefore minimize $\sum_{c\in\mathcal C}\dot W_c$ in steady state and its daily integral in transient operation:
\begin{equation}
W(\bm\theta;\bm\xi)=\int_0^T\sum_{c\in\mathcal C}\dot W_c\bigl(\alpha_c(t;\bm\theta_c),f_c(t;\bm\theta,\bm\xi)\bigr)\,dt.
\label{eq:daily_energy}
\end{equation}
Both cases require $p_{\min}\le p_j\le p_{\max}$, $1\le\alpha_c\le\bar\alpha_c$, and $0\le f_c\le\bar f_c$, where $p_{\min}$ and $p_{\max}$ are pressure limits and $\bar\alpha_c$ and $\bar f_c$ compressor ratings. These bounds hold throughout the day. We state the optimization problems in Sections~\ref{sec:steady} and~\ref{sec:transient}; the latter specifies the block and hourly Fourier schedules, admissible sets, and repeated-day initialization.

\subsection{Pressure margins}
\label{sec:margins}

The minimum-pressure requirement is local: sufficient pressure at one node cannot compensate for a shortfall at another. Subtracting the minimum pressure $p_{\min}$ from each nodal pressure gives its local margin, and the smallest local margin determines whether every node meets the limit. In steady state, we define this margin as
\begin{equation}
m(\bm\alpha;\bm\xi)=\min_{j\in\mathcal N}\bigl[p_j(\bm\alpha;\bm\xi)-p_{\min}\bigr].
\label{eq:steady_margin}
\end{equation}
The requirement holds exactly when $m\ge0$; the magnitude of a negative margin gives the worst shortfall. In transient operation, the requirement must hold throughout the day, and later recovery does not remove an earlier violation, so we take the minimum over time:
\begin{equation}
M(\bm\theta;\bm\xi)=\min_{t\in[0,T]}\min_{j\in\mathcal N_0}\bigl[p_j(t;\bm\theta,\bm\xi)-p_{\min}\bigr].
\label{eq:daily_margin}
\end{equation}
The event $M\ge0$ represents a whole day without a nodal minimum-pressure violation; its probability concerns the whole trajectory, not separate pointwise-in-time probabilities~\cite{gottlich2020}. Including node $s$, whose prescribed slack pressure lies above $p_{\min}$, leaves the minimum-pressure violation event unchanged, and node-specific limits $p_{\min,j}$ enter the same way through local margins. The maximum-pressure requirement is represented analogously by
\begin{equation}
M^{\max}(\bm\theta;\bm\xi)=\min_{t\in[0,T]}\min_{j\in\mathcal N}\bigl[p_{\max}-p_j(t;\bm\theta,\bm\xi)\bigr],
\label{eq:ceiling_margin}
\end{equation}
and holds when $M^{\max}\ge0$.

We build surrogates of the margins directly because each requirement depends on an extremum, not its location. Taking an extremum can introduce kinks when changes in compressor settings shift the limiting pressure to another node or time, even if individual pressure responses are smooth. This informs the GP choice in Section~\ref{sec:method}.
\section{Learned-margin optimization}
\label{sec:method}

\subsection{Warped Gaussian process regression}
\label{sec:warped_gp}

A GP models an unknown response by assuming that its values at any finite set of inputs are jointly Gaussian, with a mean function $\mu_0$ and covariance kernel $k$ \cite{rasmussen2006}. Conditioning on observations gives predictive distributions. For margin $y$, either $m$ or $M$, the input $\bm x$ combines decisions and prescribed conditions: $(\bm\alpha,\bm\xi)$ in steady state or $(\bm\theta,\bm\xi)$ in daily operation. Including $\bm\xi$ lets one surrogate serve different days. A finite-smoothness Mat\'ern-$5/2$ kernel with a length scale per input approximates the kinks from changes in the minimizing node or time \cite{rasmussen2006}.

A GP with a Gaussian likelihood and fixed model parameters gives Gaussian, symmetric predictions. \emph{Output warping} relaxes this restriction in physical units by modeling a monotone transformation of the response \cite{snelson2003}. To reduce skewness while accommodating margins on either side of the pressure limit, we set $g(\bm x)=\psi(y(\bm x))$ using the Yeo--Johnson transformation \cite{yeo2000}:
\begin{equation}
 \psi(y)=\begin{cases}
 [(1+y)^{\vartheta}-1]/\vartheta, & y\ge0,\\[2pt]
 -[(1-y)^{2-\vartheta}-1]/(2-\vartheta), & y<0,
 \end{cases}
\label{eq:warp}
\end{equation}
with logarithmic limits at $\vartheta=0$ and $2$ in the respective branches. We express $y$ numerically in MPa and fit $\vartheta$ before training the GP. The essential property of the Yeo--Johnson transformation is that it changes the shape of predictive uncertainty without changing which margins satisfy the pressure limit, since strict monotonicity and $\psi(0)=0$ preserve the event $y\ge0$.

For training data $\mathcal D=\{(\bm x_i,y_i)\}_{i=1}^n$, let $\bm g=[\psi(y_1),\ldots,\psi(y_n)]^\top$. With constant prior mean $\mu_0$, the posterior moments in transformed-margin units are:
\begin{equation}
\begin{aligned}
\mu_g(\bm x)&=\mu_0+\bm k_{\bm x}^{\top}\bm K^{-1}(\bm g-\mu_0\bm1),\\
\sigma_g^2(\bm x)&=k(\bm x,\bm x)-\bm k_{\bm x}^{\top}\bm K^{-1}\bm k_{\bm x}.
\end{aligned}
\label{eq:gp_posterior}
\end{equation}
Here $\bm X$ contains the training inputs, $\bm k_{\bm x}=k(\bm X,\bm x)$, and $\bm K=k(\bm X,\bm X)+\sigma_n^2I$ with fixed numerical regularization $\sigma_n^2$. Undoing training standardization puts \eqref{eq:gp_posterior} on the $g=\psi(y)$ scale; $\sigma_g$ is its posterior standard deviation. Full treatments of GP regression and output warping are given in \cite{rasmussen2006,snelson2003}.

\subsection{Epistemic chance constraints}
\label{sec:chance_constraint}

We impose the pressure requirement through a \emph{chance constraint}, using both posterior moments rather than the mean alone. Specifically, we require posterior probability at least $1-\varepsilon$ that the margin is nonnegative, with prescribed violation tolerance $\varepsilon\in(0,1)$ \cite{mesbah2016}. The GP probability $\mathbb P_{\mathcal D}$ describes epistemic uncertainty in the learned response at fixed $\bm x$, not uncertain future withdrawals.

\noindent\emph{Remark.} Since \eqref{eq:warp} preserves the zero-crossing event, warped-GP quantiles \cite{snelson2003} and the Gaussian deterministic equivalent \cite{nemirovski2006}, also used for an unwarped GP in \cite{divito2026}, give the exact equivalence under the fitted posterior
\begin{equation}
\begin{aligned}
\mathbb P_{\mathcal D}\{y(\bm x)\ge0\}&\ge1-\varepsilon\\
\Longleftrightarrow\quad\mu_g(\bm x)-z_{1-\varepsilon}\sigma_g(\bm x)&\ge0,
\end{aligned}
\label{eq:chance_equivalent}
\end{equation}
where $z_{1-\varepsilon}=\Phi^{-1}(1-\varepsilon)$ and $\Phi$ is the standard normal distribution function.

The expression $\mu_g-z_{1-\varepsilon}\sigma_g$ is the lower $\varepsilon$-quantile of the transformed margin: for $\sigma_g>0$, the posterior assigns probability $\varepsilon$ below this value. Keeping this value nonnegative limits the probability of a negative margin, rather than checking only the mean. At the studied tolerances, $z_{1-\varepsilon}>0$, so $\mu_g$ must cover the uncertainty allowance $z_{1-\varepsilon}\sigma_g$, which grows with posterior width or a smaller tolerance. Since $y$ is a minimum, one constraint covers all nodes and, for $M$, the whole day. The ceiling margin $M^{\max}$ receives a separate constraint at the same tolerance.

\subsection{Surrogate assessment and enrichment}
\label{sec:assessment}

A model may track the physical margin while attaching predictive intervals too narrow to reflect its errors. Conversely, wider intervals can improve coverage without improving the prediction. We therefore assess \emph{accuracy}, \emph{bulk coverage}, and \emph{tail coverage} together \cite{kuleshov2018,yadav2026}. We measure accuracy using root-mean-square error and the coefficient of determination $R^2$ in physical units. For a warped GP, we estimate the predictive mean by averaging inverse-transformed samples rather than inverse-transforming $\mu_g$, because a nonlinear inverse need not preserve expectations \cite{snelson2003}. Bulk coverage measures how often the central $90\%$ predictive interval contains the computed margin. Tail coverage applies the same check to the lowest decile of observed margins, where errors most affect the pressure requirement.

Passing these checks does not establish the calibration defined in Section~\ref{sec:introduction}. To examine agreement between probabilities and observed frequencies \cite{gneiting2007}, we locate each margin in its predictive distribution using the standardized error $Z_i$ and probability integral transform (PIT) $u_i$:
\begin{equation}
 Z_i=\frac{\psi(y_i)-\mu_g(\bm x_i)}{\sigma_g(\bm x_i)},\qquad u_i=\Phi(Z_i),
\label{eq:std_error}
\end{equation}
where $y_i$ is the margin computed from the reference physical model at $\bm x_i$, while $\mu_g(\bm x_i)$ and $\sigma_g(\bm x_i)>0$ are the corresponding GP mean and standard deviation on the transformed scale. The value $u_i$ is the predicted probability below the observed margin. Uniform PIT is the calibration reference; small $u_i$ indicates a lower-tail observation. When \eqref{eq:chance_equivalent} binds, violations satisfy $Z_i<-z_{1-\varepsilon}$, or $u_i<\varepsilon$, testing the constraint's tail rather than low-margin interval coverage.

In dispatch, we minimize compression work subject to \eqref{eq:chance_equivalent} and other operating limits (Sections~\ref{sec:steady}--\ref{sec:transient}). We audit each accepted decision with the reference solver or simulator and compare its margin with the GP prediction that selected it. Comparing accuracy, coverage, and standardized errors on random inputs and selected decisions tests whether reliability transfers into optimization.

The audits also supply observations precisely where the learned constraint guided the decision. We use them for \emph{enrichment}, appending audited responses from new operating days to the data for the same GP model. By \emph{conditioning}, we mean recomputing its posterior mean and covariance from the enlarged data, not refitting the model. The warp, kernel parameters, prior mean, regularization, and preprocessing remain fixed \cite{rasmussen2006}. We then freeze the surrogate and dispatch procedure and audit independent operating days. We assess minimum-pressure risk control by requiring the one-sided $95\%$ Clopper--Pearson upper bound on the violation rate among accepted dispatches not to exceed $\varepsilon$ \cite{clopper1934}. 

We also complement enrichment with empirical threshold calibration in Section~\ref{sec:transient} to bring audited violation rates closer to the declared tolerance after reoptimization.
\section{Steady-state case study}
\label{sec:steady}

\subsection{Benchmark and dispatch}
\label{sec:steady_setup}

We use the network of Fig.~\ref{fig:network}, where an operating day $\bm\xi=(d_3,d_5,p_s)$ fixes the withdrawals and supply pressure, and a dispatch $\bm\alpha$ sets the three compressor ratios. We sample the inputs independently and uniformly, with each withdrawal between $75$ and $100\%$ of its nominal value ($150$ and $125$~kg/s), and supply pressure within $\pm5\%$ of nominal (approximately $3.4474$~MPa). The pressure limits are $3$ and $6$~MPa, the maximum ratios are $(1.40,1.35,1.40)$, and the flow ratings are $(275,260,140)$~kg/s. With $\bm\Pi$ and $\bm f$ collecting squared pressures and edge flows, the \emph{physical reference problem} is
\begin{subequations}\label{eq:steady_reference}
\begin{align}
\min_{\bm\alpha,\bm\Pi,\bm f}\quad&\textstyle\sum_{c\in\mathcal C}\dot W_c(\alpha_c,f_c),\label{eq:ref_objective}\\
\text{s.t.}\quad&\eqref{eq:steady_physics},\label{eq:reference_physics}\\
&p_{\min}^2\le\Pi_j\le p_{\max}^2,\quad j\in\mathcal N,\label{eq:reference_pressure_bounds}\\
&1\le\alpha_c\le\bar\alpha_c,\quad 0\le f_c\le\bar f_c,\quad c\in\mathcal C.\label{eq:reference_bounds}
\end{align}
\end{subequations}
\begingroup
\predisplaypenalty=0
The \emph{learned problem} retains the physics, replacing only the minimum-pressure bounds by the GP constraint:
\begin{subequations}\label{eq:steady_dispatch}
\begin{align}
\min_{\bm\alpha,\bm\Pi,\bm f}\quad&\textstyle\sum_{c\in\mathcal C}\dot W_c(\alpha_c,f_c),\\
\text{s.t.}\quad&\eqref{eq:steady_physics},\ \eqref{eq:reference_bounds},\quad \Pi_j\le p_{\max}^2,\ j\in\mathcal N,\\
&\mu_g(\bm\alpha,\bm\xi)-z_{1-\varepsilon}\sigma_g(\bm\alpha,\bm\xi)\ge0.\label{eq:steady_cc_row}
\end{align}
\end{subequations}
\endgroup

We fit the base GP to $722$ physically solved Latin-hypercube inputs and assess it on $4{,}490$ held-out points. We solve the learned problem \eqref{eq:steady_dispatch} with multistart Ipopt \cite{wachter2006} at $\varepsilon\in\{1,5,10,20\}\%$ and audit accepted dispatches with the steady solver \cite{srinivasan2023}. The physical reference \eqref{eq:steady_reference} is solved separately to benchmark compression work. Audited margins from new operating days update the GP during enrichment. We evaluate the base and conditioned GPs without further updates on $500$ days, disjoint from the $500$ comparison days used for alternative designs, update rules, and formulations.

Since the margin is learned from full-physics solutions, we ask whether the pipe pressure-drop equalities can also be removed from \eqref{eq:steady_dispatch}. Those equalities also determine the flows entering the compression-work objective. Nodal balance alone fixes these flows on a tree, that is, a connected network without cycles. For the network considered here, however, balance fixes total transfer from node $2$ to node $4$, but not its split between the direct route through $P_4$ and the route through $C_2$, $P_2$, and $P_3$ (Fig.~\ref{fig:network}). Without pipe equalities, the optimizer can set $C_2$'s modeled flow to zero, eliminating its reported work at unchanged ratios. We therefore learn that physical flow with a second GP, remove \eqref{eq:pipe_law}, and replace the objective by
\begin{equation}\label{eq:reduced_objective}
\min_{\bm\alpha,\bm\Pi,\bm f}\quad
\sum_{c\in\mathcal C\setminus\{2\}}\dot W_c(\alpha_c,f_c)
+\dot W_2\bigl(\alpha_2,\mu_{f_2}(\bm\alpha,\bm\xi)\bigr),
\end{equation}
where $\mu_{f_2}$ is the GP mean used to compute work, while $f_2$ remains a balance variable. With the twelve-round margin GP, the reduction gives a negligible median paired work difference on common accepted comparison days, measured by \eqref{eq:steady_premium} with \eqref{eq:steady_dispatch} as reference. On those same days, it produces at most two additional minimum-pressure violations per tolerance compared with \eqref{eq:steady_dispatch}. This agreement, on a benchmark with an inactive physical ceiling, supports replacing the missing flow response in the objective. The transient study extends that reasoning to the objective and constraint responses supplied by the removed dynamics. With inexpensive steady-state physics, we retain the pipe relations in subsequent tests to isolate the learned pressure constraint.

\subsection{Validation at selected decisions}
\label{sec:steady_certification}

By the held-out assessment of Section~\ref{sec:assessment}, the base surrogate appears suitable for dispatch: it is accurate, and both coverages lie inside the acceptance band $[0.85,0.95]$ (Table~\ref{tab:certification}). That performance does not carry over to dispatch. At $\varepsilon=1\%$, the optimizer accepts all $500$ evaluation days, yet the physical audit finds minimum-pressure violations on $71.0\%$ of them. The intervals that covered $87.0\%$ of held-out margins cover only $12.8\%$ at these decisions.

\begin{figure}[t]
\centering
\includegraphics[width=\columnwidth]{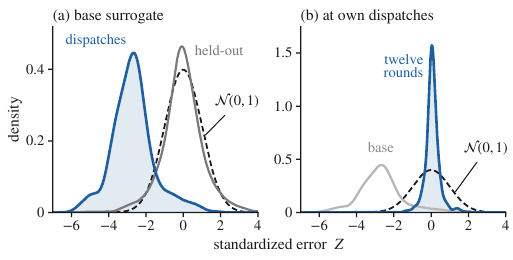}
\caption{Standardized errors at $\varepsilon=1\%$. (a) The base GP on the $4{,}490$ held-out points and its $500$ evaluation-day dispatches. (b) Base and twelve-round GPs at their own dispatches. The dashed density is standard normal.}
\label{fig:selection}
\end{figure}

A closer diagnostic reveals that the GP systematically overestimates the margin at selected decisions. Figure~\ref{fig:selection}(a) shows the standardized-error mean shifting from $-0.03$ on held-out inputs to $-2.70$ at dispatches, while its standard deviation stays similar, $1.14$ and $1.13$. The change is a shift toward overprediction, not greater scatter.

We also ask whether a fixed design could anticipate where the optimizer will operate. At the same $722$-point budget, we refit on a design concentrated near the pressure limit and on a uniform design. The former still yields $27.4\%$ violations at $\varepsilon=1\%$ on the comparison days; the latter yields $52.0\%$ despite better held-out accuracy. For these tested designs, held-out validation does not establish reliability at optimizer-selected decisions, motivating the enrichment loop below.

\begin{table}[t]
\caption{Steady GP diagnostics: central $90\%$ coverage, with tail coverage on the lowest margin decile. Decision metrics use each GP's $500$ evaluation-day dispatches at $\varepsilon=1\%$.}
\label{tab:certification}
\centering\footnotesize
\setlength{\tabcolsep}{2.6pt}
\begin{tabular}{@{}lcccccc@{}}
\toprule
&\multicolumn{4}{c}{Held-out points}&\multicolumn{2}{c}{Own dispatches}\\
\cmidrule(lr){2-5}\cmidrule(l){6-7}
GP & Cov. & Tail & $R^2$ & RMSE (kPa) & Cov. & RMSE (kPa)\\
\midrule
Base & 0.870 & 0.918 & 0.991 & 77.2 & 0.128 & 115.5\\
Four rounds & 0.858 & 0.920 & 0.992 & 71.0 & 0.988 & 5.8\\
Twelve rounds & 0.856 & 0.920 & 0.992 & 70.7 & 0.996 & 1.6\\
\bottomrule
\end{tabular}
\end{table}

\subsection{Enrichment by exact conditioning}
\label{sec:steady_enrichment}

The audits supply what these fixed designs lack: margins at the decisions the surrogate selects. Each round dispatches $32$ new days at the four tolerances and conditions on their audited margins. Four and twelve rounds add $512$ and $1{,}536$ observations and meet the minimum-pressure risk criterion at all four tolerances (Table~\ref{tab:steady_rates}). After examining all twelve rounds, we selected round four as the earliest to meet the criterion on both day sets. The reported rates reuse those selection days, so they are not an independent test of a budget fixed in advance. The results show that further conditioning improves decision-point accuracy while held-out diagnostics barely move (Table~\ref{tab:certification}). It also makes the policy more conservative at $5\%$ and above, with errors narrower than the assumed normal law (Fig.~\ref{fig:selection}(b)).

To price risk control, let $\dot W_i^{\mathrm{GP}}$ and $\dot W_i^{\mathrm{ref}}$ be the audited total compression-work rates of the learned policy and physical reference \eqref{eq:steady_reference} on common accepted days $i\in\mathcal I$. Their \emph{median paired premium} (\%) is
\begin{equation}\label{eq:steady_premium}
\operatorname{median}_{i\in\mathcal I}
\left[100\,\frac{\dot W_i^{\mathrm{GP}}-\dot W_i^{\mathrm{ref}}}
{\dot W_i^{\mathrm{ref}}}\right].
\end{equation}
At $\varepsilon=1\%$, the premium falls from $0.70\%$ after four rounds to $0.33\%$ after twelve despite greater conservatism.

The update itself matters, not merely the additional data. On one $2{,}258$-point design from a parallel enrichment stream, we compare conditioning, refitting the kernel and mean with the warp held, and full refitting. At $\varepsilon=1\%$ on the $500$ comparison days, violations are $0.2\%$, $16.8\%$, and $20.2\%$, respectively. Conditioning gives fewer violations than either refit at every studied tolerance. Also, conditioning took about $0.2$~s for one factorization, against roughly $50$~min for a recorded refit under shared workstation load. Neither refit meets the risk criterion at any tolerance, whereas conditioning meets it at every studied tolerance.

\begin{table}[t]
\caption{Audited minimum-pressure violation rates (\%) on the same $500$ evaluation days, with one-sided $95\%$ Clopper--Pearson upper bounds in parentheses.}
\label{tab:steady_rates}
\centering\footnotesize
\setlength{\tabcolsep}{3.2pt}
\begin{tabular}{@{}lrrrr@{}}
\toprule
&\multicolumn{4}{c}{Declared tolerance $\varepsilon$}\\
GP & $1\%$ & $5\%$ & $10\%$ & $20\%$\\
\midrule
Base & 71.0 (74.3) & 87.8 (90.1) & 91.8 (93.7) & 96.2 (97.5)\\
Four rounds & 0.0 (0.6) & 0.6 (1.5) & 2.0 (3.4) & 5.8 (7.8)\\
Twelve rounds & 0.0 (0.6) & 0.0 (0.6) & 0.2 (0.9) & 1.8 (3.1)\\
\bottomrule
\end{tabular}
\end{table}
\section{Transient case study}
\label{sec:transient}

\subsection{Daily inputs and compressor schedules}
\label{sec:transient_controls}
Daily operation couples compressor decisions through linepack, so a schedule must be assessed through the pressure trajectory it produces. We use the network and operating limits of Section~\ref{sec:steady_setup}, now with time-varying withdrawals over $T=24$~h. For withdrawal node $j\in\{3,5\}$, we set
\begin{equation}
 d_j(t)=\bar d_j\nu_j\bigl[1-\beta_j\cos\omega(t-t_j)\bigr],
 \qquad \omega=\frac{2\pi}{T},
\label{eq:withdrawal_profile}
\end{equation}
where $\bar d_j$ is its nominal withdrawal, $\nu_j\in[0.75,1.00]$ sets the relative daily mean, $\beta_j\in[0.08,0.26]$ the relative swing, and $t_j\in[-3,3]$~h the time of minimum withdrawal. These parameters for nodes $3$ and $5$, together with constant supply pressure $p_s\in[3.275,3.620]$~MPa, give $\bm\xi=(\nu_3,t_3,\beta_3,\nu_5,t_5,\beta_5,p_s)$. Demand follows five-minute linear interpolation of \eqref{eq:withdrawal_profile} within control intervals.

\begin{figure}[!ht]
\centering
\includegraphics[width=\columnwidth]{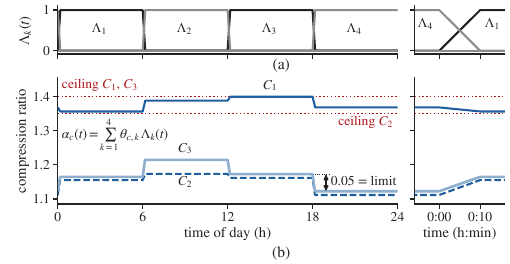}
\caption{Block schedules. (a) Periodic basis functions and the midnight transition. (b) Applied ratios for the example in Fig.~\ref{fig:transient_day}; dotted lines mark ratio limits. The inset resolves the ten-minute ramp.}
\label{fig:block_family}
\end{figure}

To keep the GP input dimension manageable, we parameterize schedules with a limited set of basis coefficients~\cite{diehl2024noc}, independent of the finer time grid used to resolve the dynamics. We study two families: the \emph{block family} and the \emph{hourly Fourier family}. The block family uses four six-hour levels per compressor, joined by ten-minute ramps:
\begin{equation}
 \alpha_c(t;\bm\theta_c)=\sum_{k=1}^{4}\theta_{c,k}\Lambda_k(t).
\label{eq:block_family}
\end{equation}
Each periodic trapezoid $\Lambda_k$ rises from zero to one in the first ten minutes of block $k$, stays at one for the remainder, and falls in the first ten minutes of the next block. The basis functions sum to one, making $\theta_{c,k}$ the block level. The day begins with the ramp from the fourth level to the first (Fig.~\ref{fig:block_family}). Three compressors require twelve coefficients.

\begin{figure}[!ht]
\centering
\includegraphics[width=\columnwidth]{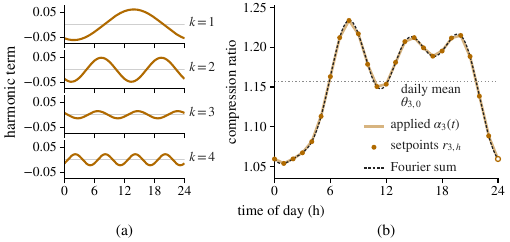}
\caption{Hourly Fourier schedules. (a) The four harmonic contributions for $C_3$. (b) Adding the daily mean gives the Fourier sum (dashed); hourly samples (dots) are joined linearly to produce the applied ratio (solid).}
\label{fig:fourier_family}
\end{figure}

The hourly Fourier family instead generates hourly setpoints from a daily mean and four cosine--sine pairs:
\begin{equation}
\begin{aligned}
 r_{c,h}&=\theta_{c,0}+\sum_{k=1}^{4}\left[\theta_{c,2k-1}\cos\frac{2\pi kh}{24}
                  +\theta_{c,2k}\sin\frac{2\pi kh}{24}\right],\\
 \alpha_c(t;\bm\theta_c)&=\sum_{h=0}^{24}r_{c,h}H_h(t),\qquad r_{c,24}=r_{c,0}.
\end{aligned}
\label{eq:hourly_family}
\end{equation}
Here $h=0,\ldots,23$ indexes the hourly setpoints. One mean coefficient $\theta_{c,0}$ and eight harmonic coefficients determine all $24$ setpoints for each compressor, giving nine coefficients per compressor and $27$ in total. Repeating $r_{c,0}$ at hour $24$ closes the daily schedule. The hat function $H_h$ equals one at hour $h$, zero at the other hourly knots, and varies linearly between adjacent knots. The resulting schedule joins the setpoints by straight lines rather than applying the continuous Fourier sum (Fig.~\ref{fig:fourier_family}). Thus, optimizing the coefficients selects the compressor schedule within this family.

For every compressor, we impose ratio bounds and cyclic increments through
\begin{equation}
\begin{aligned}
\mathcal A_{\mathrm{blk}}=\{\bm\theta:\;&1\le\theta_{c,k}\le\bar\alpha_c,\quad
 |\theta_{c,k+1}-\theta_{c,k}|\le\Delta\alpha\},\\
\mathcal A_{\mathrm{hr}}=\{\bm\theta:\;&1\le r_{c,h}\le\bar\alpha_c,\quad
 |r_{c,h+1}-r_{c,h}|\le\Delta\alpha\},
\end{aligned}
\label{eq:transient_admissible}
\end{equation}
with $\Delta\alpha=0.05$, $\theta_{c,5}=\theta_{c,1}$, and $r_{c,24}=r_{c,0}$. These constraints define bounded convex polytopes in coefficient space, and interpolation preserves the ratio bounds. The basis and ramp restrictions make the families nonnested.

\begin{figure*}[!t]
\centering
\includegraphics[width=\textwidth]{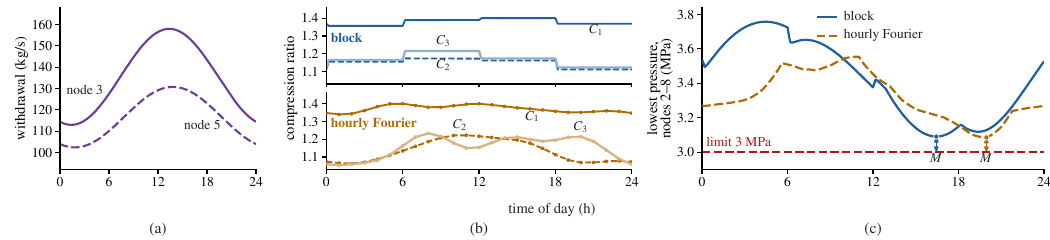}
\caption{One operating day under the twelve-round maps at $\varepsilon=5\%$, using nominal thresholds. (a) Withdrawals. (b) Applied compressor ratios. (c) Lowest nodal pressures from the two-second audits, with whole-day margins $M$. This day illustrates the families, not their aggregate performance.}
\label{fig:transient_day}
\end{figure*}

\subsection{Learned dispatch and physical evaluation}
\label{sec:transient_dispatch}
We use GasTranSim as the reference physical model for each day and schedule, integrating \eqref{eq:gas_dynamics} with the explicit staggered-grid method~\cite{gyrya2019}. After a twelve-hour pre-run at the day's initial boundary values, we repeat each day and schedule from the preceding final state until nodal pressures satisfy $\max_{j\in\mathcal N}|p_j(T)-p_j(0)|\le10^{-4}$~MPa or six passes have run. We use the last pass, evaluating margins at every two-second step for both families; one-second fresh-day re-audits assess numerical sensitivity.

The steady pipe-removal comparison in Section~\ref{sec:steady_setup} supplied the missing objective response with a second GP. Here we learn all responses needed for dispatch, so the optimizer evaluates schedules without integrating or retaining pressure and flow trajectories. Alongside the margins $M$ and $M^{\max}$ in \eqref{eq:daily_margin}--\eqref{eq:ceiling_margin}, we therefore learn daily work $W$ from \eqref{eq:daily_energy} and the flow extrema $f_c^{\mathrm{pk}}=\max_{t\in[0,T]}f_c(t)$ and $f_c^{\mathrm{lo}}=\min_{t\in[0,T]}f_c(t)$. The pressure margins use the warped GPs of Section~\ref{sec:warped_gp}; work and flow maps are unwarped. Writing $\mu_g^M,\sigma_g^M$ and $\mu_g^{M^{\max}},\sigma_g^{M^{\max}}$ for transformed-margin moments, we solve
\begin{subequations}\label{eq:transient_program}
\begin{align}
\min_{\bm\theta\in\mathcal A}\quad&\mu_W(\bm\theta,\bm\xi),\label{eq:transient_objective}\\
\text{s.t.}\quad&\mu_g^M(\bm\theta,\bm\xi)-z_{1-\varepsilon}\sigma_g^M(\bm\theta,\bm\xi)\ge0,\label{eq:transient_floor}\\
&\mu_g^{M^{\max}}(\bm\theta,\bm\xi)-z_{1-\varepsilon}\sigma_g^{M^{\max}}(\bm\theta,\bm\xi)\ge0,\label{eq:transient_ceiling}\\
&\mu_{f_c^{\mathrm{pk}}}(\bm\theta,\bm\xi)\le\bar f_c,\quad
 \mu_{f_c^{\mathrm{lo}}}(\bm\theta,\bm\xi)\ge0,\quad c\in\mathcal C,\label{eq:transient_flows}
\end{align}
\end{subequations}
where $\mathcal A$ is the corresponding set in \eqref{eq:transient_admissible}. Each pressure limit receives its own chance constraint; the flow limits use posterior means alone. Consistent with the transient formulation in~\cite{kazi2024acc}, we include compression work in the objective but impose no separate power limits. Thus, the optimizer sees a finite list of responses rather than a trajectory, while each physical audit still checks the whole day (Fig.~\ref{fig:transient_day}).

Combining schedule coefficients with the seven day parameters gives $19$ GP inputs for the block family and $34$ for the hourly family. Excluding failed pre-runs leaves $1{,}199$ block and $2{,}046$ hourly minimum-margin targets, including imputed failures. The ceiling-margin, daily-work, and flow-extrema GPs instead use the $908$ completed block and $1{,}530$ completed hourly simulations. Held-out assessment uses $305$ completed block and $1{,}511$ completed hourly cases.

We solve \eqref{eq:transient_program} by multistart sequential quadratic programming with analytic GP derivatives. Each round offers $32$ new days at $\varepsilon\in\{1,5,10,20\}\%$ and conditions on accepted-schedule audits with fitted quantities fixed. We compare the base and round maps without updates on $500$ days per family, drawn independently and uniformly over the day-parameter ranges. We report violations among accepted schedules and compare audited work only on days both policies dispatch.

\subsection{Block schedules}
\label{sec:transient_block}
The base block-margin GP is accurate on completed held-out cases ($R^2=0.976$), with bulk coverage $0.944$ but tail coverage $0.774$, below the acceptance band. At selected schedules, however, its central interval covers $96$--$98\%$ of audited margins. The standardized errors have root-mean-square magnitude $0.67$--$0.76$, below the unit reference for calibration. Thus, despite a modest optimistic bias, predictive distributions are wide relative to the errors. Unlike the steady case, this conservatism limits the initial violations (Table~\ref{tab:block_results}).

\begin{table}[!ht]
\centering
\caption{Block minimum-pressure violation rates (\%) among accepted schedules, with one-sided 95\% upper bounds.}
\label{tab:block_results}
\setlength{\tabcolsep}{3.3pt}
\small
\begin{tabular}{@{}lrrrr@{}}
\toprule
 &\multicolumn{4}{c}{Declared tolerance $\varepsilon$}\\
Maps&$1\%$&$5\%$&$10\%$&$20\%$\\
\midrule
Base&0.2 (1.1)&2.4 (4.0)&7.7 (10.1)&20.5 (23.9)\\
Four rounds&0.0 (0.7)&0.2 (1.1)&2.4 (4.0)&7.8 (10.3)\\
Twelve rounds&0.0 (0.8)&0.8 (1.9)&2.0 (3.6)&8.1 (10.7)\\
\bottomrule
\end{tabular}
\end{table}

After four and twelve rounds of conditioning, the policies meet the minimum-pressure risk criterion at all four tolerances on these days, while their selected-decision intervals remain conservative. Enrichment also reduces median paired work relative to the base by $1.7$--$5.4\%$ after four rounds and $2.5$--$7.5\%$ after twelve, using the daily-energy counterpart of \eqref{eq:steady_premium}. Also, acceptance, the fraction of offered days receiving a dispatch, falls from $90$--$91\%$ to $79$--$81\%$ after twelve rounds. Every newly declined day had exceeded a compressor flow rating in the base policy's audit. Conditioning therefore improves the decisions and work predictions. We note, however, that the work saving applies to days both policies dispatch, not to serving every offered day.

\subsection{Hourly schedules and threshold calibration}
\label{sec:transient_hourly}
The hourly margin GP gives similarly strong held-out accuracy ($R^2=0.971$), with bulk and tail coverages $0.931$ and $0.829$. Conditioning substantially reduces violations at selected schedules, but neither four nor twelve rounds meets the risk criterion at any tolerance (Table~\ref{tab:hourly_results}, first three rows). After twelve rounds, standardized errors have nearly unit spread (standard deviation $0.91$--$0.95$) but retain an optimistic mean offset of $-0.33$ to $-0.50$. The numerically active minimum-pressure constraint connects these violations to the remaining lower-tail mismatch through \eqref{eq:std_error}.

\begin{table}[!ht]
\centering
\caption{Hourly minimum-pressure violation rates (\%) among accepted schedules, with one-sided 95\% upper bounds. Empirical: 1,000 fresh days, 1~s audit; other rows: 500 evaluation days, 2~s audit.}
\label{tab:hourly_results}
\setlength{\tabcolsep}{1.9pt}
\small
\begin{tabular}{@{}lrrrr@{}}
\toprule
 &\multicolumn{4}{c}{Declared tolerance $\varepsilon$}\\
Policy&$1\%$&$5\%$&$10\%$&$20\%$\\
\midrule
Base&15.0 (18.0)&36.8 (40.6)&56.3 (60.1)&74.4 (77.6)\\
Four rounds&0.8 (1.9)&5.8 (7.9)&14.5 (17.4)&28.7 (32.2)\\
Twelve rounds&2.3 (3.7)&8.6 (11.0)&17.0 (20.0)&31.8 (35.5)\\
\midrule
Empirical&1.03 (1.73)&4.71 (5.98)&9.36 (11.03)&16.45 (18.51)\\
\bottomrule
\end{tabular}
\end{table}

We therefore adapt \emph{split-conformal quantile selection}~\cite{angelopoulos2023} to set the allowance from audited errors. For each tolerance, let $S_i=-Z_i$ for the $N_\varepsilon$ accepted calibration decisions under the twelve-round maps and nominal Gaussian threshold $z_{1-\varepsilon}$. Larger scores indicate greater overestimation relative to predictive width. With $S_{(k)}$ the $k$th smallest score, we set
\begin{equation}
\begin{aligned}
 k_\varepsilon&=\lceil(N_\varepsilon+1)(1-\varepsilon)\rceil,
 \qquad q_\varepsilon=S_{(k_\varepsilon)},\\
 \mu_g^M(\bm\theta,\bm\xi)&-q_\varepsilon\sigma_g^M(\bm\theta,\bm\xi)\ge0.
\end{aligned}
\label{eq:empirical_threshold}
\end{equation}
Only the minimum-pressure threshold changes, not the posterior. The $1{,}500$ calibration days also informed design and map selection, so reuse and reoptimization preclude the ordinary split-conformal guarantee~\cite{fannjiang2022}. With maps and thresholds frozen, we reoptimize on $1{,}000$ fresh days. Table~\ref{tab:hourly_results} adds the empirical policy's one-second audits to the three $500$-day, two-second evaluations.

The thresholds increase the nominal allowances by factors of $1.23$--$1.50$. One-second fresh-day audits show closer alignment across all four tolerances, including a reduction from $10.17\%$ under the nominal threshold on the same fresh days to $4.71\%$ at $5\%$. Paired reductions are significant at every tolerance (largest one-sided McNemar $p=1.9\times10^{-6}$). Median paired work increases by $0.82$--$0.99\%$ against the nominal rule, measured at two seconds.
\section{Conclusion and Future Work}
\label{sec:conclusion}
We develop a learned-margin formulation for compressor scheduling, coupling a scalar pressure chance constraint with assessment at selected decisions. By modeling the whole-day minimum directly and representing schedules through block levels or Fourier-generated hourly setpoints, we optimize daily operation without embedding pressure and flow trajectories. Our steady-state diagnostic shows that held-out accuracy and coverage do not establish dispatch reliability, motivating loop-based enrichment by exact conditioning. Conditioning meets the risk criterion at every studied tolerance, unlike matched-data refits, with about four orders of magnitude less update time in the shared-workstation comparison. In transient operation, it recovers conservative block risk control but leaves the hourly Fourier policy optimistic. This residual hourly optimism motivates empirical threshold calibration, bringing fresh-day rates closer to all declared tolerances for less than $1\%$ additional paired work.

Our current study, however, has limitations. Power limits are omitted, while learned flow and separate pressure-ceiling constraints still admit violations, so joint operational feasibility is not established. Conservative risk control also differs from probabilistic consistency, motivating calibration that avoids unnecessary compression while respecting the declared tolerance. Fixed block and truncated Fourier bases restrict schedule flexibility, while finer control discretizations can make the GP input dimension prohibitive. Without a matching full-physics transient optimizer, we do not establish a computational speedup, but our formulation removes the full space-time state discretization from the optimization problem. We thus view adaptive GPs and calibration for probabilistically consistent dispatch, adaptive harmonic selection for flexibility, and Bayesian graph neural networks for scalable, interpretable prediction as future extensions.

\section*{ACKNOWLEDGMENT}
T.~A.~Adebiyi performed this work as an Applied Machine Learning Fellow and Graduate Research Assistant at Los Alamos National Laboratory. Research presented in this article was supported by the Laboratory Directed Research and Development program of Los Alamos National Laboratory under project number 3W500A-YWR0-0000. Research conducted at Los Alamos National Laboratory is done under the auspices of the National Nuclear Security Administration of the U.S. Department of Energy under Contract No.~89233218CNA000001.

\bibliographystyle{IEEEtran}
\bibliography{references}

\begin{thebibliography}{10}
\providecommand{\url}[1]{#1}
\csname url@samestyle\endcsname
\providecommand{\newblock}{\relax}
\providecommand{\bibinfo}[2]{#2}
\providecommand{\BIBentrySTDinterwordspacing}{\spaceskip=0pt\relax}
\providecommand{\BIBentryALTinterwordstretchfactor}{4}
\providecommand{\BIBentryALTinterwordspacing}{\spaceskip=\fontdimen2\font plus
\BIBentryALTinterwordstretchfactor\fontdimen3\font minus \fontdimen4\font\relax}
\providecommand{\BIBforeignlanguage}[2]{{%
\expandafter\ifx\csname l@#1\endcsname\relax
\typeout{** WARNING: IEEEtran.bst: No hyphenation pattern has been}%
\typeout{** loaded for the language `#1'. Using the pattern for}%
\typeout{** the default language instead.}%
\else
\language=\csname l@#1\endcsname
\fi
#2}}
\providecommand{\BIBdecl}{\relax}
\BIBdecl

\bibitem{srinivasan2023}
S.~Srinivasan, K.~Sundar, V.~Gyrya, and A.~Zlotnik, ``Numerical solution of the steady-state network flow equations for a nonideal gas,'' \emph{IEEE Trans. Control Netw. Syst.}, vol.~10, no.~3, pp. 1449--1461, 2023.

\bibitem{zlotnik2015}
A.~Zlotnik, M.~Chertkov, and S.~Backhaus, ``Optimal control of transient flow in natural gas networks,'' in \emph{Proc. 54th IEEE Conf. Decis. Control}, 2015, pp. 4563--4570.

\bibitem{arrigo2022}
A.~Arrigo, M.~Dol{\'a}nyi, K.~Bruninx, and J.-F. Toubeau, ``Machine learning for improved gas network models in coordinated energy systems,'' 2022, arXiv:2209.12731.

\bibitem{wang2025}
S.~Wang, W.~Wu, C.~Lin, Q.~Wang, S.~Xu, and B.~Chen, ``Physics-informed recurrent network for state-space modeling of gas pipeline networks,'' 2025, arXiv:2502.07230.

\bibitem{roald2020}
L.~A. Roald, K.~Sundar, A.~Zlotnik, S.~Misra, and G.~Andersson, ``An uncertainty management framework for integrated gas-electric energy systems,'' \emph{Proc. IEEE}, vol. 108, no.~9, pp. 1518--1540, 2020.

\bibitem{dvorkin2022}
V.~Dvorkin, A.~Ratha, P.~Pinson, and J.~Kazempour, ``Stochastic control and pricing for natural gas networks,'' \emph{IEEE Trans. Control Netw. Syst.}, vol.~9, no.~1, pp. 450--462, 2022.

\bibitem{kazi2024sfv}
S.~R. Kazi, S.~Misra, S.~Tokareva, K.~Sundar, and A.~Zlotnik, ``Stochastic finite volume method for uncertainty management in gas pipeline network flows,'' in \emph{Proc. 63rd IEEE Conf. Decis. Control}, 2024, pp. 82--89.

\bibitem{jiang2026}
D.~Jiang, J.~Garcke, O.~Akca, J.~Hollnagel, B.~Klaassen, M.~Anvari, and J.~M{\"u}ller-Kirchenbauer, ``Accelerating gas-network feasibility screening with a physics-informed graph neural network surrogate,'' 2026, arXiv:2607.13610.

\bibitem{dubourg2011}
V.~Dubourg, B.~Sudret, and J.-M. Bourinet, ``Reliability-based design optimization using kriging surrogates and subset simulation,'' \emph{Struct. Multidiscip. Optim.}, vol.~44, no.~5, pp. 673--690, 2011.

\bibitem{rasmussen2006}
C.~E. Rasmussen and C.~K.~I. Williams, \emph{Gaussian Processes for Machine Learning}.\hskip 1em plus 0.5em minus 0.4em\relax Cambridge, MA: MIT Press, 2006.

\bibitem{mesbah2016}
A.~Mesbah, ``Stochastic model predictive control: An overview and perspectives for future research,'' \emph{IEEE Control Syst. Mag.}, vol.~36, no.~6, pp. 30--44, 2016.

\bibitem{wiebe2022}
J.~Wiebe, I.~Cec{\'i}lio, J.~Dunlop, and R.~Misener, ``A robust approach to warped {Gaussian} process-constrained optimization,'' \emph{Math. Program.}, vol. 196, pp. 805--839, 2022.

\bibitem{divito2026}
V.~{Di Vito}, K.~Sundar, F.~Fioretto, and D.~Deka, ``Stability-constrained {AC} optimal power flow--{A Gaussian} process-based approach,'' \emph{Sustain. Energy Grids Netw.}, vol.~46, p. 102300, 2026.

\bibitem{snelson2003}
E.~Snelson, C.~E. Rasmussen, and Z.~Ghahramani, ``Warped {Gaussian} processes,'' in \emph{Adv. Neural Inf. Process. Syst.}, vol.~16, 2003.

\bibitem{nemirovski2006}
A.~Nemirovski and A.~Shapiro, ``Convex approximations of chance constrained programs,'' \emph{SIAM J. Optim.}, vol.~17, no.~4, pp. 969--996, 2006.

\bibitem{gneiting2007}
T.~Gneiting, F.~Balabdaoui, and A.~E. Raftery, ``Probabilistic forecasts, calibration and sharpness,'' \emph{J. Roy. Statist. Soc. B}, vol.~69, no.~2, pp. 243--268, 2007.

\bibitem{kuleshov2018}
V.~Kuleshov, N.~Fenner, and S.~Ermon, ``Accurate uncertainties for deep learning using calibrated regression,'' in \emph{Proc. 35th Int. Conf. Mach. Learn.}, vol.~80, 2018, pp. 2796--2804.

\bibitem{yadav2026}
A.~Yadav, T.~A. Adebiyi, and R.~Zhang, ``Calibrating scientific foundation models with inference-time stochastic attention,'' in \emph{Adv. Neural Inf. Process. Syst.}, 2026, to be published, arXiv:2604.19530.

\bibitem{allen2025}
S.~Allen, J.~Koh, J.~Segers, and J.~Ziegel, ``Tail calibration of probabilistic forecasts,'' \emph{J. Amer. Statist. Assoc.}, vol. 120, no. 552, pp. 2796--2808, 2025.

\bibitem{fannjiang2022}
C.~Fannjiang, S.~Bates, A.~N. Angelopoulos, J.~Listgarten, and M.~I. Jordan, ``Conformal prediction under feedback covariate shift for biomolecular design,'' \emph{Proc. Nat. Acad. Sci.}, vol. 119, no.~43, p. e2204569119, 2022.

\bibitem{gyrya2019}
V.~Gyrya and A.~Zlotnik, ``An explicit staggered-grid method for numerical simulation of large-scale natural gas pipeline networks,'' \emph{Appl. Math. Model.}, vol.~65, pp. 34--51, 2019.

\bibitem{gottlich2020}
S.~G{\"o}ttlich, O.~Kolb, and K.~Lux, ``Chance-constrained optimal inflow control in hyperbolic supply systems with uncertain demand,'' \emph{Optim. Control Appl. Methods}, vol.~42, no.~2, pp. 566--589, 2021.

\bibitem{yeo2000}
I.-K. Yeo and R.~A. Johnson, ``A new family of power transformations to improve normality or symmetry,'' \emph{Biometrika}, vol.~87, no.~4, pp. 954--959, 2000.

\bibitem{clopper1934}
C.~J. Clopper and E.~S. Pearson, ``The use of confidence or fiducial limits illustrated in the case of the binomial,'' \emph{Biometrika}, vol.~26, no.~4, pp. 404--413, 1934.

\bibitem{wachter2006}
A.~W{\"a}chter and L.~T. Biegler, ``On the implementation of an interior-point filter line-search algorithm for large-scale nonlinear programming,'' \emph{Math. Program.}, vol. 106, no.~1, pp. 25--57, 2006.

\bibitem{diehl2024noc}
\BIBentryALTinterwordspacing
M.~Diehl and S.~Gros, ``{Numerical Optimal Control},'' Draft book, 11 Oct. 2024. [Online]. Available: \url{https://www.syscop.de/files/2024ws/NOC/book-NOCSE.pdf}
\BIBentrySTDinterwordspacing

\bibitem{kazi2024acc}
S.~R. Kazi, K.~Sundar, and A.~Zlotnik, ``Dynamic optimization and optimal control of hydrogen blending operations in natural gas networks,'' in \emph{Proc. Amer. Control Conf.}, 2024, pp. 5357--5363.

\bibitem{angelopoulos2023}
A.~N. Angelopoulos and S.~Bates, ``Conformal prediction: A gentle introduction,'' \emph{Found. Trends Mach. Learn.}, vol.~16, no.~4, pp. 494--591, 2023.

\end{thebibliography}
\end{document}